\documentclass[11pt]{article}

\usepackage[latin1]{inputenc}\usepackage{epsfig}\usepackage{amsfonts}

\title{
  {\huge Fourier Theory on the Complex Plane IV} \\
  Representability of Real Functions \\
  by their Fourier Coefficients }

\author{
  \Large Jorge L. deLyra \\
  Department of Mathematical Physics \\
  Physics Institute \\
  University of São Paulo }

\date{May 4, 2015}

\newcommand{\ii}{\mbox{\boldmath$\imath$}}

\newcommand{\at}[2]{\left.\rule{0em}{3ex}\right[_{\,#1}^{\,#2}}

\newcommand{\e}[1]{\,{\rm e}^{#1}}

\begin{document}\maketitle

\begin{abstract}
  \noindent
  The results presented in this paper are refinements of some results
  presented in a previous paper. Three such refined results are presented.
  The first one relaxes one of the basic hypotheses assumed in the
  previous paper, and thus extends the results obtained there to a wider
  class of real functions. The other two relate to a closer examination of
  the issue of the representability of real functions by their Fourier
  coefficients. As was shown in the previous paper, in many cases one can
  recover the real function from its Fourier coefficients even if the
  corresponding Fourier series diverges almost everywhere. In such cases
  we say that the real function is still representable by its Fourier
  coefficients. Here we establish a very weak condition on the Fourier
  coefficients that ensures the representability of the function by those
  coefficients. In addition to this, we show that any real function that
  is absolutely integrable can be recovered almost everywhere from, and
  hence is representable by, its Fourier coefficients, regardless of
  whether or not its Fourier series converges. Interestingly, this also
  provides proof for a conjecture proposed in the previous paper.
\end{abstract}

\section{Introduction}

In a previous paper~\cite{FTotCPI} we have developed a correspondence
between, on one side, the Fourier series and Fourier coefficients of real
functions on the interval $[-\pi,\pi]$ and, on the other side, a complex
analytic structure within the open unit disk, consisting of a set of inner
analytic functions and their complex Taylor series. The reader is referred
to that paper for the definition of many of the concepts and notations
used in this one. In many cases this correspondence leads to the
formulation of expressions involving modified trigonometric series that
can converge very fast to a given real function, even when the Fourier
series of that function diverges or converges very slowly, as shown
in~\cite{FTotCPII}.

In order to establish the correspondence described above,
in~\cite{FTotCPI} we assumed that the Definite Parity (DP) real functions
$f(\theta)$ under examination are such that their Fourier series and the
corresponding Fourier Conjugate (FC) series converge together on at least
one point on the interval $[-\pi,\pi]$, a domain which, in the context of
the correspondence to be established, is mapped onto the unit circle of
the complex plane. In this paper we will show that one can relax that
hypothesis, replacing it by a much weaker one.

Following the development in~\cite{FTotCPI}, we will deal here only with
real functions that have definite parity properties, which we will call
Definite Parity real functions or DP real functions for short. Since any
real function $f(\theta)$ in the interval $[-\pi,\pi]$, without any
restriction, can be separated into even and odd parts,

\noindent
\begin{eqnarray*}
  f(\theta)
  & = &
  f_{\rm c}(\theta)+f_{\rm s}(\theta),
  \\
  f_{\rm c}(\theta)
  & = &
  \frac{f(\theta)+f(-\theta)}{2},
  \\
  f_{\rm s}(\theta)
  & = &
  \frac{f(\theta)-f(-\theta)}{2},
\end{eqnarray*}

\noindent
where

\noindent
\begin{eqnarray*}
  f_{\rm c}(\theta)
  & = &
  f_{\rm c}(-\theta),
  \\
  f_{\rm s}(\theta)
  & = &
  -f_{\rm s}(-\theta),
\end{eqnarray*}

\noindent
we can restrict the discussion to Definite Parity (DP) real functions
without any loss of generality. This condition implies that the
corresponding inner analytic function $w(z)$ has the property that it
reduces to a real function on the interval $(-1,1)$ of the real axis, as
discussed in~\cite{FTotCPI}. For simplicity, we will also assume that
$f(\theta)$ is a zero-average function, since adding a constant function
to $f(\theta)$ is a trivial operation that does not significantly affect
the issues under discussion here. This condition implies that the
corresponding inner analytic function $w(z)$ has the property that
$w(0)=0$, as discussed in~\cite{FTotCPI}.

In Subsection~$8.1$ of~\cite{FTotCPI} we have shown that, if there is any
singularity of an analytic function $w(z)$ inside the open unit disk, then
the sequence of Taylor-Fourier coefficients of its Taylor series diverges
to infinity exponentially fast on the unit circle, as a function of the
series index. In this paper we will discuss the converse of this
statement. We will show here that the mere absence of an
exponentially-fast divergence of the sequence of Fourier coefficients on
the unit circle is enough to ensure the convergence of the corresponding
complex power series inside the open unit disk, thus leading to the
definition of a corresponding inner analytic function.

A note about the concept of integrability of real functions is in order at
this point. What we mean by integrability in this paper is integrability
in the sense of Lebesgue, with the use of the usual Lebesgue measure. We
will assume that all the real functions at issue here are measurable in
the Lebesgue measure. Therefore whenever we speak of real functions, it
should be understood that we mean Lebesgue-measurable real functions. We
will use the following result from the theory of measure and integration:
for real functions defined within a compact interval, which are Lebesgue
measurable, integrability and absolute integrability are equivalent
conditions~\cite{RealAnalysis}. Therefore we will use the concepts of
integrability and of absolute integrability interchangeably.

\section{Weakening the Convergence Hypothesis}

Let $f(\theta)$ be a DP real function that has zero average value, with
$\theta\in[-\pi,\pi]$. We will assume that $f(\theta)$ is an integrable
function, so that its Fourier coefficients $\alpha_{k}$ and $\beta_{k}$
exist, since they are given by

\noindent
\begin{eqnarray*}
  \alpha_{k}
  & = &
  \frac{1}{\pi}
  \int_{-\pi}^{\pi}d\theta\,
  f(\theta)
  \cos(k\theta),
  \\
  \beta_{k}
  & = &
  \frac{1}{\pi}
  \int_{-\pi}^{\pi}d\theta\,
  f(\theta)
  \sin(k\theta).
\end{eqnarray*}

\noindent
Note that $\alpha_{0}=0$ because $f(\theta)$ is assumed to be a
zero-average function. We will name such coefficients $a_{k}$, for
$k=1,2,3,\ldots,\infty$, irrespective of whether $f(\theta)$ is even or
odd, and thus of whether it gives origin respectively to a cosine or sine
series. Let us assume that the $a_{k}$ coefficients satisfy the following
hypothesis: given a sequence $a_{k}$ of coefficients, there exists a
positive real function $g(k)$ with the property that

\begin{equation}\label{ratecondition1}
  \lim_{k\to\infty}g(k)
  =
  1,
\end{equation}

\noindent
and such that for all $k\geq 1$

\begin{equation}\label{ratecondition2}
  \left|\frac{a_{k+1}}{a_{k}}\right|
  \leq
  g(k).
\end{equation}

\noindent
What this means is that we assume that the ratio of the absolute values of
two consecutive coefficients is bounded from above by a function that
tends to one in the $k\to\infty$ limit. Note that the $k\to\infty$ limit
of the ratio of coefficients itself may not even exist. Note also that, if
this condition holds for a given DP real function, then it automatically
holds for the corresponding FC real function as well, since both have the
same coefficients and the condition is imposed only on these coefficients.
Following the development described in~\cite{FTotCPI}, given $f(\theta)$
with such properties we may now construct the two DP trigonometric series

\noindent
\begin{eqnarray*}
  S_{\rm c}
  & = &
  \sum_{k=1}^{\infty}
  a_{k}\cos(k\theta),
  \\
  S_{\rm s}
  & = &
  \sum_{k=1}^{\infty}
  a_{k}\sin(k\theta),
\end{eqnarray*}

\noindent
which are the FC series of one another, and then the complex power series

\noindent
\begin{eqnarray*}
  S_{z}
  & = &
  \sum_{k=1}^{\infty}
  a_{k}\rho^{k}
  \left[
    \cos(k\theta)+\ii\sin(k\theta)
  \right]
  \\
  & = &
  \sum_{k=1}^{\infty}
  a_{k}z^{k},
\end{eqnarray*}

\noindent
where $z=\rho\exp(\ii\theta)$, with $\rho\geq 0$, so that $S_{\rm c}$ and
$S_{\rm s}$ are respectively the real and imaginary parts of $S_{z}$ for
$\rho=1$. Note that the condition expressed by
Equations~(\ref{ratecondition1}) and~(\ref{ratecondition2}) does not, by
itself, imply the convergence of the series $S_{\rm c}$ and $S_{\rm s}$ on
the unit circle, since sequences of coefficients $a_{k}$ that do not go to
zero as $k\to\infty$ may satisfy it. Therefore, many sequences of
coefficients leading to Fourier series that diverge on the unit circle
satisfy that condition. If we now use the ratio criterion to analyze the
convergence of the power series $S_{z}$ inside the open unit disk, we get

\begin{displaymath}
  \left|\frac{a_{k+1}z^{k+1}}{a_{k}z^{k}}\right|
  =
  \left|\frac{a_{k+1}}{a_{k}}\right|\rho.
\end{displaymath}

\noindent
Our hypothesis about the ratio of the $a_{k}$ coefficients now leads to
the relation, inside the open unit disk,

\begin{displaymath}
  \left|\frac{a_{k+1}z^{k+1}}{a_{k}z^{k}}\right|
  \leq
  \rho g(k).
\end{displaymath}

\noindent
Now, given any value of $\rho<1$, since the limit of $g(k)$ for
$k\to\infty$ is one, it follows that the same limit of $\rho g(k)$ is
strictly less than one. Therefore, we may conclude that there is a value
$k_{m}$ of $k$ such that, if $k>k_{m}$, then $\rho g(k)<1$, so that the
ratio is less that one, thus implying that we have

\begin{displaymath}
  \left|\frac{a_{k+1}z^{k+1}}{a_{k}z^{k}}\right|
  <
  1,
\end{displaymath}

\noindent
for $k>k_{m}$. This implies that the ratio criterion is satisfied and
therefore that the series $S_{z}$ converges at such points. Since this is
true for any $\rho<1$, we may conclude that the open unit disk is the
maximum disk of convergence of the power series $S_{z}$. Therefore, the
series converges to a complex analytic function $w(z)$ there, which is an
inner analytic function according to the definition given
in~\cite{FTotCPI}.

The results established in~\cite{FTotCPI} then imply that $f(\theta)$ can
now be obtained almost everywhere over the unit circle as the $\rho\to 1$
limit of the real or imaginary part of $w(z)$, as the case may be. Note
that the series $S_{z}$, and hence the series $S_{\rm c}$ and $S_{\rm s}$,
may still be divergent over the whole unit circle. This does not affect
the recovery of the real function from its Fourier coefficients in the
manner just described. We have therefore shown that the correspondence
established in~\cite{FTotCPI} holds without the hypothesis that $S_{\rm
  c}$ and $S_{\rm s}$ be convergent together on at least at one point of
the unit circle, so long as the Taylor-Fourier coefficients $a_{k}$
satisfy the weaker condition expressed by Equations~(\ref{ratecondition1})
and~(\ref{ratecondition2}).

\section{Representability by the Fourier Coefficients}

Again we start with an arbitrary integrable DP real function $f(\theta)$
that has zero average value, with $\theta\in[-\pi,\pi]$. Let us now assume
that this function is such that the corresponding Fourier coefficients
$a_{k}$ satisfy the condition that

\begin{equation}\label{expolimit1}
  \lim_{k\to\infty}
  |a_{k}|\e{-Ck}
  =
  0,
\end{equation}

\noindent
for all real $C>0$. What this means is that $a_{k}$ may or may not go to
zero as $k\to\infty$, may approach a non-zero real number, and may even
diverge to infinity as $k\to\infty$, so long as it does not do so
exponentially fast. This includes therefore not only the sequences of
Fourier coefficients corresponding to all possible convergent Fourier
series, but many sequences that correspond to Fourier series that diverge
almost everywhere. In fact, it even includes sequences of coefficients
that cannot be obtained at all from a real function, such as the sequence
$a_{k}=1/\pi$ for all $k\geq 1$, which is associated to the Dirac delta
``function'' $\delta(\theta)$, as shown in~\cite{FTotCPI} and as will be
discussed in Section~\ref{discexts} of this paper. It is therefore a very
weak condition indeed.

Before anything else, let us establish a preliminary result, namely that
the condition in Equation~(\ref{expolimit1}) implies that we also have

\begin{equation}\label{expolimit2}
  \lim_{k\to\infty}
  |a_{k}|k^{p}\e{-Ck}
  =
  0,
\end{equation}

\noindent
for all real powers $p>0$. This is just a formalization of the well-known
fact that the negative-exponent real exponential function goes to zero
faster than any positive power goes to infinity, as $k\to\infty$. We may
write the function of $k$ on the left-hand side as

\begin{displaymath}
  |a_{k}|k^{p}\e{-Ck}
  =
  |a_{k}|\e{p\ln(k)}\e{-Ck}.
\end{displaymath}

\noindent
Recalling the properties of the logarithm, we now observe that, given an
arbitrary real number $A>0$, there is always a sufficiently large value
$k_{m}$ of $k$ above which $\ln(k)<Ak$. A simple proof can be found in
Appendix~\ref{appproplog}. Due to this we may write, for all $k>k_{m}$,

\begin{displaymath}
  |a_{k}|k^{p}\e{-Ck}
  <
  |a_{k}|\e{pAk}\e{-Ck},
\end{displaymath}

\noindent
since the exponential is a monotonically increasing function. If we choose
$A=C/(2p)$, which is positive and not zero, we get that, for all
$k>k_{m}$,

\noindent
\begin{eqnarray*}
  |a_{k}|k^{p}\e{-Ck}
  & < &
  |a_{k}|\e{Ck/2}\e{-Ck}
  \\
  & = &
  |a_{k}|\e{-Ck/2}.
\end{eqnarray*}

\noindent
According to our hypothesis about the coefficients $a_{k}$, the
$k\to\infty$ limit of the expression in the right-hand side is zero for
any strictly positive value of $C'=C/2$, so that taking the $k\to\infty$
limit we establish our preliminary result,

\begin{displaymath}
  \lim_{k\to\infty}
  |a_{k}|k^{p}\e{-Ck}
  =
  0,
\end{displaymath}

\noindent
for all real $C>0$ and all real $p>0$. If we now construct the complex
power series $S_{z}$ as before, using the coefficients $a_{k}$, we are in
a position to show that it is absolutely convergent inside the open unit
disk. In order to do this we consider the real power series
$\overline{S}_{z}$ of the absolute values of the terms of the series
$S_{z}$, which we write as

\noindent
\begin{eqnarray*}
  \overline{S}_{z}
  & = &
  \sum_{k=1}^{\infty}
  |a_{k}|\rho^{k}
  \\
  & = &
  \sum_{k=1}^{\infty}
  |a_{k}|\e{k\ln(\rho)}.
\end{eqnarray*}

\noindent
Since $\rho<1$ inside the open unit disk, the logarithm shown is strictly
negative, and we may put $\ln(\rho)=-C$ with real $C>0$. We can now see
that, according to our hypothesis about the coefficients $a_{k}$, the
terms of this series go to zero exponentially fast as $k\to\infty$. This
already suffices to establish its convergence, but we may easily make this
more explicit, writing

\noindent
\begin{eqnarray*}
  \overline{S}_{z}
  & = &
  \sum_{k=1}^{\infty}
  |a_{k}|\e{-Ck}
  \\
  & = &
  \sum_{k=1}^{\infty}
  \frac{k^{2}|a_{k}|\e{-Ck}}{k^{2}}.
\end{eqnarray*}

\noindent
According to our preliminary result in Equation~(\ref{expolimit2}) the
numerator shown goes to zero as $k\to\infty$, and therefore above a
sufficiently large value $k_{m}$ of $k$ it is less that one, so we may
write that

\noindent
\begin{eqnarray*}
  \overline{S}_{z}
  & = &
  \sum_{k=1}^{k_{m}}
  |a_{k}|\e{-Ck}
  +
  \sum_{k=k_{m}+1}^{\infty}
  \frac{k^{2}|a_{k}|\e{-Ck}}{k^{2}}
  \\
  & < &
  \sum_{k=1}^{k_{m}}
  |a_{k}|\e{-Ck}
  +
  \sum_{k=k_{m}+1}^{\infty}
  \frac{1}{k^{2}}.
\end{eqnarray*}

\noindent
The first term on the right-hand side is a finite sum and therefore is
finite, and the second term can be bounded from above by a convergent
asymptotic integral on $k$, so that we have

\noindent
\begin{eqnarray*}
  \overline{S}_{z}
  & < &
  \sum_{k=1}^{k_{m}}
  |a_{k}|\e{-Ck}
  +
  \int_{k_{m}}^{\infty}dk\,
  \frac{1}{k^{2}}
  \\
  & = &
  \sum_{k=1}^{k_{m}}
  |a_{k}|\e{-Ck}
  +
  \frac{-1}{k}\at{k_{m}}{\infty}
  \\
  & = &
  \sum_{k=1}^{k_{m}}
  |a_{k}|\e{-Ck}
  +
  \frac{1}{k_{m}}.
\end{eqnarray*}

\noindent
It follows that $\overline{S}_{z}$, which is a real sum of positive terms,
so that its partial sums form a monotonically increasing sequence, is
bounded from above and is therefore convergent. It then follows that
$S_{z}$ is absolutely convergent and therefore convergent. Since this is
valid for all $\rho<1$, we may conclude that the maximum disk of
convergence of $S_{z}$ is the open unit disk. We may now recover
$f(\theta)$ as the $\rho\to 1$ limit of the real or imaginary part of
$w(z)$, as the case may be, almost everywhere on the unit circle, as was
mentioned before and shown in~\cite{FTotCPI}.

This provides therefore a very general condition on the Fourier
coefficients of the real functions that ensures that the correspondence
established in~\cite{FTotCPI} holds. This then ensures that the real
functions can be recovered from, and therefore can be represented by,
their Fourier coefficients. Note that the condition in
Equation~(\ref{expolimit1}) can be considered as an even weaker form of
the condition discussed in the previous section. We are therefore ready to
state our first important conclusion:

\vspace{3ex}

\noindent
\parbox{\textwidth} {\bf\boldmath If the sequence of Fourier coefficients
  $a_{k}$ of a DP real function $f(\theta)$ satisfies the weak condition
  in Equation~(\ref{expolimit1}), so that the coefficients do not diverge
  to infinity exponentially fast or faster when $k\to\infty$, then the
  complex power series constructed from them converges to an inner
  analytic function inside the open unit disk, from which the real
  function can be recovered almost everywhere.}

\vspace{3ex}

\noindent
Note that, although we formulated this condition in terms of the Fourier
coefficients $a_{k}$ of a given DP real function, the fact that $a_{k}$
are the Fourier coefficients of the function has in fact not been used at
all. Therefore, the conclusion is valid for any sequence of coefficients
that satisfies Equation~(\ref{expolimit1}), regardless of whether or not
they can be obtained as the Fourier coefficients of some real function.

\section{Analytic Criterion for Real Functions}

Finally, let us establish a simple analytic condition over the real
functions that ensures that they are representable by their Fourier
coefficients. If we assume that $f(\theta)$ is absolutely integrable, so
that the integral

\begin{displaymath}
  \frac{1}{\pi}
  \int_{-\pi}^{\pi}d\theta\,
  |f(\theta)|
  =
  F
\end{displaymath}

\noindent
exists and is a finite real number $F$, then it follows that we have for
the Fourier coefficients, taking as an example the case of the cosine
series,

\noindent
\begin{eqnarray*}
  |a_{k}|
  & = &
  \left|
    \frac{1}{\pi}
    \int_{-\pi}^{\pi}d\theta\,
    f(\theta)
    \cos(k\theta)
  \right|
  \\
  & \leq &
  \frac{1}{\pi}
  \int_{-\pi}^{\pi}d\theta\,
  |f(\theta)|
  |\cos(k\theta)|
  \\
  & \leq &
  \frac{1}{\pi}
  \int_{-\pi}^{\pi}d\theta\,
  |f(\theta)|
  \\
  & = &
  F,
\end{eqnarray*}

\noindent
where we used the triangular inequalities. It follows that we have, for
all $k\geq 1$,

\begin{displaymath}
  |a_{k}|
  \leq
  F.
\end{displaymath}

\noindent
Since we thus see that the Fourier coefficients of $f(\theta)$ are bounded
within the interval $[-F,F]$, for all $k\geq 1$, it follows that they
cannot diverge to infinity as $k\to\infty$, and therefore that they
satisfy our hypothesis in Equation~(\ref{expolimit1}), namely that

\begin{displaymath}
  \lim_{k\to\infty}
  |a_{k}|\e{-Ck}
  =
  0,
\end{displaymath}

\noindent
for all real $C>0$. The same result can be established in a similar way
for the case of the sine series, of course. It therefore follows that
$f(\theta)$ is representable by its Fourier coefficients.

If we go back to a more general function $f(\theta)$ that has both even
and odd parts, since the result holds for both parts, since we may also
add a constant term without changing the result, and since we have limited
ourselves to Lebesgue-measurable real functions within a compact interval,
for which integrability and absolute integrability are one and the same
concept, we are ready to state our second important conclusion:

\vspace{3ex}

\noindent
\parbox{\textwidth} {\bf\boldmath Any integrable real function $f(\theta)$
  defined within $[-\pi,\pi]$ is representable by its Fourier
  coefficients, and can be recovered from them almost everywhere in that
  interval, regardless of whether or not the corresponding Fourier series
  converges.}

\vspace{3ex}

\noindent
It is an interesting observation that this provides an answer to the
conjecture proposed in~\cite{FTotCPI}, about whether or not there are any
integrable real functions such that their sequences of Fourier
coefficients $a_{k}$ give rise to complex power series $S_{z}$ which are
strongly divergent, that is, that have at least one singular point
strictly within the open unit disk. The answer, according to the proof
worked out here, is that there are none, as expected.

\section{Discussion and Extensions}\label{discexts}

In their use of Fourier series for the solution of physical problems,
physicists often, and quite successfully, simply ignore convergence issues
of the Fourier series involved. With just a bit of common sense and a
willingness to accept approximate results, they just plow ahead with their
calculations, and that seldom leads them into serious trouble. This is
true even if one includes the occasional appearance in these calculations
of the Fourier series of singular objects such as Dirac's delta
``function'', since there is usually not much difficulty in interpreting
the divergences in physical terms.

The results established here can be viewed as an explanation of this
rather remarkable fact. Functions that are useful in physics applications
are always at the very least Lebesgue-measurable, and most often at least
locally integrable almost everywhere. Extremely pathological functions are
not of any use in such circumstances. Therefore, for all real functions of
interest in physics applications, there is in effect an underlying
analytic structure that firmly anchors all the operations which are
performed on the Fourier series, mapping them onto corresponding and much
safer operations on the inner-analytic functions within the open unit
disk. These operations may include anything from the basic arithmetic
operations to integration and differentiation, and so on.

The representation of real functions by their Fourier coefficients, which
in physics usually goes by the name of ``representation in momentum
space'', often can be interpreted directly in terms of physical concepts.
In fact, this alternative representation is frequently found to be the
more important and fundamental one. In terms of the mathematical structure
that we are dealing with here, it is clear that this state of affairs
relates closely to the fact that the analytic structure within the open
unit disk is quite clearly the more fundamental aspect of this whole
mathematical structure. In the usual physics parlance, that analytic
structure is an exact and universal ``regulator'' for all integrable real
functions.

It may be possible to further extend the result presented in the previous
section. Observe that if the real function $f(\theta)$ is limited and
integrable, then its Fourier coefficients $a_{k}$ are also limited, and
thus it is obvious that they satisfy the condition in
Equation~(\ref{expolimit1}), so that they do not diverge exponentially
fast with $k$. However, the function does not have to be limited in order
for the coefficients to satisfy that condition. The function may diverge
to infinity at an isolated point, so long as the asymptotic integral
around that point exists. If the function diverges to infinity at a point
$\theta_{0}$ in such a way that the integrals

\noindent
\begin{eqnarray*}
  I_{\ominus}
  & = &
  \int_{\theta_{\ominus}}^{\theta_{0}}d\theta\,
  f(\theta),
  \\
  I_{\oplus}
  & = &
  \int_{\theta_{0}}^{\theta_{\oplus}}d\theta\,
  f(\theta),
\end{eqnarray*}

\noindent
exist and are finite for some $\theta_{\ominus}$ and $\theta_{\oplus}$,
with $\theta_{\ominus}<\theta_{0}<\theta_{\oplus}$, and where there are no
other hard singularities of $f(\theta)$ within the intervals
$[\theta_{\ominus},\theta_{0})$ and $(\theta_{0},\theta_{\oplus}]$, then
the function is integrable and thus representable by its Fourier
coefficients. One may also have a finite number of such isolated points
without disturbing these properties. One may even consider the inclusion
of a denumerable infinity of such points, so long as the numerical series
resulting from the sum of the contributions of all the singular points to
the integral is absolutely convergent, since otherwise the function cannot
be considered to be integrable.

Another interesting extension of the results presented in the two previous
sections would be one leading of the inclusion in the structure of
singular objects which are not real functions, such as extended functions
or distributions in the sense of Schwartz~\cite{DistTheory}, represented
by their distribution kernels, such as Dirac's delta ``function''. This
extension seems to be relatively straightforward, and these objects are
routinely dealt with, without too much trouble, in physics
applications. These singular objects are also associated to inner analytic
functions.

Note that, since the Fourier coefficients of any absolutely integrable
real function are necessarily limited, these sequences of coefficients
form only a small subset of the set of all the sequences of coefficients
that satisfy the condition stated in Equation~(\ref{expolimit1}), since
many sequences of coefficients which are not limited may satisfy that
condition. This shows once again that the condition on the coefficients is
in fact very weak. Besides, that condition may be applied to any sequence
of real coefficients, whether or not they are the Fourier coefficients of
a real function. Therefore, the condition includes much more than just
real functions, since there may be many sequences of coefficients $a_{k}$
that satisfy it but that are not obtainable from a real function on the
unit circle, as its sequence of Fourier coefficients.

In other words, there are inner analytic functions within the open unit
disk that correspond to definite sequences of coefficients $a_{k}$
satisfying our hypothesis, but that are not related to a real function on
the unit circle. These inner analytic functions have at least one hard
singularity over the unit circle, such as a simple pole, or harder. A hard
singularity is defined in~\cite{FTotCPII} as a point at which the limit of
the inner-analytic function does not exist, or diverges to infinity. A
related gradation in terms of degrees of hardness is also defined there.
One important example of this is the inner analytic function associated to
the Dirac delta ``function'', which was given in~\cite{FTotCPI}. Given a
point $z_{1}$ on the unit circle, the very simple analytic function

\begin{displaymath}
  w_{\delta}(z)
  =
  \frac{1}{2\pi}
  -
  \frac{1}{\pi}\,
  \frac{z}{z-z_{1}},
\end{displaymath}

\noindent
which has a simple pole on the unit circle, is an extended inner analytic
function, that is, an inner analytic function rotated by the angle
$\theta_{1}$ associated to $z_{1}=\exp(\ii\theta_{1})$ and with the
constant shown added to it. This analytic function within the open unit
disk is such that the delta ``function'' can be obtained as the limit to
the unit circle of its real part,

\begin{displaymath}
  \delta(\theta-\theta_{1})
  =
  \lim_{\rho\to 1}
  \Re[w_{\delta}(z)],
\end{displaymath}

\noindent
as was shown in~\cite{FTotCPI} with basis on the properties that define
the delta ``function''. It would be interesting to further investigate
this extension, which would probably include every integrable object which
is not a real function. In particular, it would be interesting to
investigate whether or not there is a condition such as the condition of
absolute integrability which is a sufficient condition, and possibly also
a necessary condition, for the representability of such objects by their
Fourier coefficients.

\section{Conclusions}

Two very weak conditions were established leading to the existence of
corresponding inner analytic functions, and thus to the representability
of real functions or other objects defined on the unit circle by their
sequence of Fourier coefficients. One of them is a condition on the
sequence of coefficients, the other is an analytic condition on the real
functions. The first one, stated in the most general way possible, reads:

\vspace{3ex}

\noindent
\parbox{\textwidth} {\bf\boldmath If a sequence of real coefficients
  $a_{k}$ satisfies the weak condition in Equation~(\ref{expolimit1}), so
  that the coefficients do not diverge to infinity exponentially fast or
  faster when $k\to\infty$, then the complex power series constructed from
  them converges to an inner analytic function inside the open unit disk.}

\vspace{3ex}

\noindent
This is a very weak condition on the sequence of coefficients, leading to
the representability of the object related to it by an inner analytic
function. The object at issue may be an integrable real function, or it
may be a singular object that has an integrability concept associated to
it. The typical example of such singular objects is Dirac's delta
``function'' and the Schwartz distribution associated to it. This is
therefore a very general condition, which in fact extrapolates the strict
realm of real functions. The second one reads:

\vspace{3ex}

\noindent
\parbox{\textwidth} {\bf\boldmath Any Lebesgue-measurable integrable real
  function $f(\theta)$ defined within $[-\pi,\pi]$ is representable by its
  Fourier coefficients, and can be recovered from them almost everywhere
  in that interval, regardless of whether or not the corresponding Fourier
  series converges.}

\vspace{3ex}

\noindent
This is a very weak analytical condition on the real functions, leading to
the representability of each real function by a corresponding inner
analytic function. Observe that, so long as the real functions at issue
are Lebesgue-measurable, this is the statement that {\em every} integrable
function is representable by its sequence of Fourier coefficients. All
that is really needed is that these Fourier coefficients, and hence the
integrals that give them, exist.

If one accepts the premise that the condition of integrability is
necessary for the very existence of the sequence of Fourier coefficients,
and since it is seen here to be also sufficient for the representability
of the function by those coefficients, we may say that what we have here
is in essence a necessary and sufficient condition for the
representability of real functions by their Fourier coefficients, as
defined by the usual integrals. This is similar in nature to the as yet
open problem of establishing a necessary and sufficient condition for the
convergence of the Fourier series of a real function. If we reinterpret
this last one as a condition for the representability of the function by
its series, we see that such a necessary and sufficient condition can be
achieved if one exchanges the condition of representability by the series
for a more general one, involving the representability by the sequence of
coefficients, and the recovery of the real function from them via the
construction of an analytic function within the open unit disk. However,
note that this condition is only necessary if one assumes that the Fourier
coefficients must be given by the usual integrals over the periodic
interval.

Although the process of recovery of the real function as a limit of the
corresponding inner analytic function is, by itself, not algorithmic in
nature, it can lead to algorithmic solutions in at least some cases. The
results obtained here simply guarantee the existence of the inner analytic
function, and at least in principle the possibility of the recovery of the
real function through the limit of that inner analytic function to the
unit circle. In sufficiently simple cases, namely when the inner analytic
function has only a finite number of sufficiently soft known singularities
on the unit circle, it is possible to devise expressions involving
modified trigonometric series that converge to the function, as was shown
in~\cite{FTotCPII}. In this case, one acquires an algorithmic method for
the calculation of the real function to any desired level of precision,
and thus for the practical recovery of the real function from the Fourier
coefficients, which can be made good enough for any practical purpose.

\appendix

\section{A Property of the Logarithm}\label{appproplog}

Let us show that the logarithm has the property that given an arbitrary
real number $A>0$, there is always a sufficiently large value $k_{m}$ of
the integer $k$ above which $\ln(k)<Ak$. We simply promote $k$ to a
continuous real variable $x>0$ and consider the function

\begin{displaymath}
  h(x)
  =
  Ax-\ln(x).
\end{displaymath}

\noindent
Is is quite clear that this function diverges to positive infinity as
$x\to 0$. If we calculate the first and second derivatives of $h(x)$ we
get

\noindent
\begin{eqnarray*}
  h'(x)
  & = &
  A-\frac{1}{x},
  \\
  h''(x)
  & = &
  \frac{1}{x^{2}}.
\end{eqnarray*}

\noindent
It is now clear that there is a single critical point where the first
derivative is zero, that is where $h'(x_{0})=0$, given by $x_{0}=1/A$. At
this point we have for the second derivative

\begin{displaymath}
  h''(x_{0})
  =
  A^{2},
\end{displaymath}

\noindent
which is positive, implying that the critical point is a local minimum.
Since there is no other minimum, maximum or inflection point, it becomes
clear that the function must decrease from positive infinity as $x$
increases from zero, go through the point of minimum at $x_{0}$, and then
increase without limit as $x\to\infty$. At this point of minimum we have
for the function itself,

\begin{displaymath}
  h(x_{0})
  =
  1+\ln(A).
\end{displaymath}

\noindent
It follows that, if $h(x_{0})>0$, then we must have $h(x)>0$ for all
$x>0$. This corresponds to $A>1/e$. On the other hand, if $h(x_{0})\leq
0$, then there are two solutions $x_{1}$ and $x_{2}$ of the equation
$h(x)=0$, that coincide if $h(x_{0})=0$. This corresponds to $0<A\leq
1/e$. In this case the function is positive within the interval
$(0,x_{1})$, negative within $(x_{1},x_{2})$ and positive for $x>x_{2}$.

Therefore, for all possible values of $A$ there is a value $x_{m}$ of $x$,
either $x_{m}=0$ or $x_{m}=x_{2}$, such that for $x>x_{m}$ the function
$h(x)$ is positive, and therefore

\begin{displaymath}
  Ax
  >
  \ln(x).
\end{displaymath}

\noindent
Translating the statement back in terms of $k$, we have that given any
real number $A>0$, there is a minimum value $k_{m}$ of the integer $k$
such that $Ak>\ln(k)$. Thus the required statement is established.

\section{Proof of Uniqueness Almost Everywhere}

Given a sequence of Fourier coefficients $a_{k}$ of an integrable DP real
function $f(\theta)$ with zero average, let us show that they uniquely
characterize that real function almost everywhere, that is, up to a
zero-measure function. The point of the proof is to show that this is true
independently of the convergence of the Fourier series of the function,
hence including the cases in which that series diverges almost everywhere.
Since two real functions with opposite parities clearly must always have
different sequences of Fourier coefficients, it is enough to show that the
result holds within the two sets of real functions with the same parity.
Although the definition of the integral in itself will not be used
directly, the argument will involve the concept of a zero-measure
function, and therefore we should always think in terms of the Lebesgue
integral and of the usual Lebesgue measure, for conceptual reasons as well
as for the sake of generality. Hence, when we talk here of an integrable
function, we mean integrable in the sense of Lebesgue.

Imagine then that we are given two different integrable DP real functions
$f_{1}(\theta)$ and $f_{2}(\theta)$, both with zero average and both with
the same parity, and consider their difference

\begin{displaymath}
  g(\theta)
  =
  f_{1}(\theta)-f_{2}(\theta),
\end{displaymath}

\noindent
which of course has the same parity as $f_{1}(\theta)$ and
$f_{2}(\theta)$, and which is also a zero-average function. Let us assume
that both $f_{1}(\theta)$ and $f_{2}(\theta)$ have the same sequence of
Fourier coefficients $a_{k}$, for $k=1,2,3,\ldots,\infty$. It follows that
$g(\theta)$ has all its Fourier coefficients $b_{k}$ equal to zero, since
from the relation above we clearly have

\begin{displaymath}
  b_{k}
  =
  a_{k}-a_{k},
\end{displaymath}

\noindent
for all $k\geq 1$. This means that the function $g(\theta)$ has zero
scalar product with all the elements of the Fourier basis, and thus that
it is orthogonal to all of them. However, as shown in~\cite{FTotCPI}, that
basis satisfies a completeness relation,

\begin{equation}\label{completeness}
  \delta(\theta-\theta_{1})
  =
  \frac{1}{2\pi}
  +
  \frac{1}{\pi}
  \sum_{k=1}^{\infty}
  \left[
    \cos(k\theta)
    \cos(k\theta_{1})
    +
    \sin(k\theta)
    \sin(k\theta_{1})
  \right],
\end{equation}

\noindent
for all $\theta$ and all $\theta_{1}$, and is therefore complete. From
this it follows that $g(\theta)$, being orthogonal to all the elements of
the basis, must be a zero-measure function, that is, a function which is
zero-measure equivalent to the identically zero function. In order to show
this, we simply multiply the equation above by $g(\theta)$ and integrate
on $\theta$ over the periodic interval,

\noindent
\begin{eqnarray*}
  \int_{-\pi}^{\pi}d\theta\,
  g(\theta)
  \delta(\theta-\theta_{1})
  & = &
  \int_{-\pi}^{\pi}d\theta\,
  g(\theta)\,
  \frac{1}{2\pi}
  +
  \\
  &   &
  +
  \int_{-\pi}^{\pi}d\theta\,
  g(\theta)\,
  \frac{1}{\pi}
  \sum_{k=1}^{\infty}
  \left[
    \cos(k\theta)
    \cos(k\theta_{1})
    +
    \sin(k\theta)
    \sin(k\theta_{1})
  \right]
  \\
  & = &
  \frac{1}{2}
  \left[
    \frac{1}{\pi}
    \int_{-\pi}^{\pi}d\theta\,
    g(\theta)
  \right]
  +
  \\
  &   &
  +
  \sum_{k=1}^{\infty}
  \left\{
    \left[
      \frac{1}{\pi}
      \int_{-\pi}^{\pi}d\theta\,
      g(\theta)\,
      \cos(k\theta)
    \right]
    \cos(k\theta_{1})
  \right.
  +
  \\
  &   &
  \hspace{2.3em}
  +
  \left.
    \left[
      \frac{1}{\pi}
      \int_{-\pi}^{\pi}d\theta\,
      g(\theta)\,
      \sin(k\theta)
    \right]
    \sin(k\theta_{1})
  \right\}.
\end{eqnarray*}

\noindent
The integral on the left-hand side can be calculated through the use of
the properties of the delta ``function'', which were given and
demonstrated in~\cite{FTotCPI}. The three expressions within square
brackets on the right-hand side are the various Fourier coefficients of
the function $g(\theta)$. The first one is zero because $g(\theta)$, just
like $f_{1}(\theta)$ and $f_{2}(\theta)$, is a zero-average function. The
other two sequences of coefficients are identically zero because the
scalar products of $g(\theta)$ and the elements of the basis are all zero.
Besides, one of the sequences contains the coefficients of the part of the
basis which has the wrong parity, and thus these coefficients must all be
zero by a parity argument. Therefore the whole right-hand side is zero and
hence we have

\begin{displaymath}
  g(\theta_{1})
  =
  0,
\end{displaymath}

\noindent
which is valid for almost all values of $\theta_{1}$. Therefore
$g(\theta)$ is zero almost everywhere. In other words, $f_{1}(\theta)$ and
$f_{2}(\theta)$ differ only by a zero-measure function, and are thus equal
to each other almost everywhere.

We may conclude from this that, if $f_{1}(\theta)$ and $f_{2}(\theta)$ are
to be different by more than a zero-measure function, then their sequences
of Fourier coefficients must be different, and therefore the coefficients
uniquely characterize the originating function, up to a zero-measure
function. Observe that nothing in this argument involves the convergence
of the Fourier series of the real functions. In particular, the proof of
completeness of the Fourier basis presented in~\cite{FTotCPI} and
expressed by the relation in Equation~(\ref{completeness}) was obtained
directly from the complex analytic structure within the open unit disk,
and it also does not involve in any way the convergence of the Fourier
series on the unit circle.

\section{Zero-Measure Equivalence Classes}

The whole structure we are examining here induces one to think that it
would be a reasonable thing to do if we decided to group all the real
functions into zero-measure equivalence classes. We could consider as
equivalent, at least from the point of view of the physics applications,
two real functions which differ by a zero-measure function. For this
purpose, a zero-measure real function is a function whose absolute
integral over $[-\pi,\pi]$ is zero, or that has zero integral on any
closed sub-interval of the interval $[-\pi,\pi]$.

Then one would consider in a group all functions that are zero measure in
this sense, such as the identically null function. This class could then
be represented simply by that particular element, $f(\theta)\equiv0$,
which is quite clearly the smoothest element within that class. Given any
real function with non-zero measure, all other functions that differ from
it by a zero-measure function would be in the same equivalence class.
Instead of considering all real functions, we could formulate everything
that was discussed here in terms of these equivalence classes.

Since the Fourier coefficients are defined by integrals, it is immediately
clear that a given sequence of Fourier coefficients, if it belongs to any
real function at all, belongs to one such equivalence class, rather than
to the individual functions. Therefore, there is also some mathematical
sense to such a classification. The proposed method of representation of
the real functions, as limits to the unit circle of real or imaginary
parts of inner analytic functions within the open unit disk, gives us then
a clear definition of a preferred or central element of each class: that
particular real function which is obtained at the limit just described.
Being a limit from a perfectly smooth analytic function, this is clearly
the smoothest element of the class.

One could consider this classification process, and the representation of
each zero-measure equivalence class by its smoothest element, as a process
of elimination of what we might call the irrelevant pathologies of the
real functions. The only pathologies that would remain would be those that
have a definite effect on the integral of the function. This certainly
makes sense in terms of the physics applications, but it might be of some
intrinsic mathematical interest as well.


\begin{thebibliography}{99}

\bibitem{FTotCPI} J. L. deLyra, ``Fourier Theory on the Complex Plane I --
  Conjugate Pairs of Fourier Series and Inner Analytic Functions'', arXiv:
  1409.2582.

\bibitem{FTotCPII} J. L. deLyra, ``Fourier Theory on the Complex Plane II
  -- Weak Convergence, Classification and Factorization of
  Singularities'', arXiv: 1409.4435.

\bibitem{RealAnalysis} W. Rudin, ``Principles of Mathematical Analysis'',
  McGraw-Hill, third edition, 1976; ISBN-13: 978-0070542358, ISBN-10:
  007054235X.

  H. Royden, ``Real Analysis'', Prentice-Hall Inc., third edition, 1988;
  ISBN-13: 978-0024041517, ISBN-10: 0024041513.

\bibitem{DistTheory} L. Schwartz, ``Théorie des Distributions'', {\bf
    1-2}, Hermann, Paris, (1951).

  M. J. Lighthill, ``Introduction to Fourier Analysis and Generalized
  Functions'', Cambridge University Press, (1959). ISBN 0-521-09128-4.

\end{thebibliography}
\end{document}